\documentclass[11pt]{amsart}
\usepackage[dvips]{graphicx}
\usepackage[all]{xy}
\usepackage{amsmath}
\usepackage{amssymb}

\newtheorem{theorem}{Theorem}
\newtheorem{lemma}{Lemma}
\theoremstyle{definition}
\newtheorem{definition}[theorem]{Definition}
\newtheorem{remark}{Remark}
\numberwithin{equation}{section}

\theoremstyle{corollary}

\numberwithin{equation}{section}
\theoremstyle{example}

\newtheorem{corollary}{Corollary}
\theoremstyle{proposition}
\newtheorem{proposition}{Proposition}
\newcommand\R{\mathbb{R}}

\newcommand\C{\mathbb{C}}

\renewcommand\P{\mathbb{P}}
\newcommand\Z{\mathbb{Z}}

\newcommand\G{{\rm G}}

\newcommand\Ker{{\rm Ker}}

\newcommand\tr{{\rm tr}}

\newcommand\diag{{\rm diag}}

\newcommand\Aut{{\rm Aut}}

\newcommand\SO{{\rm SO}}
\renewcommand\O{{\rm O}}
\newcommand\Spin{{\rm Spin}}
\newcommand\SU{{\rm SU}}
\newcommand\SL{{\rm SL}}
\newcommand\Sp{{\rm Sp}}
\newcommand\GL{{\rm GL}}

\newcommand\U{{\rm U}}

\newcommand\Ea{{\rm E}}
\newcommand\Fa{{\rm F}}
\begin{document}
%
%
\title{Homogeneous Lagrangian submanifolds}
\author{Lucio Bedulli and Anna Gori}
\address{Dipartimento di Matematica - Universit\`a di Bologna\\
Piazza di Porta S. Donato 5\\40126 Bologna\\Italy} 
\email{bedulli@math.unifi.it}
\address{Dipartimento di Matematica e Appl.\ per l'Architettura - Universit\`a di Firenze\\
Piazza Ghiberti 27\\50100 Firenze\\Italy}
\email{gori@math.unifi.it}
\begin{abstract}
We characterize isometric actions on compact K\"ahler manifolds admitting a Lagrangian orbit, 
describing under which condition the Lagrangian orbit is unique.
We furthermore give the complete classification of simple groups acting on the complex
projective space with a Lagrangian orbit, and we give the explicit list  of these orbits.
\end{abstract}
\maketitle
\section*{Introduction}
%
%
%
%
     A Lagrangian submanifold of a $2n$-dimensional symplectic
 manifold $(M,\omega)$ is an $n$-dimensional submanifold on which
 the symplectic form $\omega$ vanishes. Lagrangian submanifolds
 play an important role in symplectic geometry and topology.\\
 In the K\"ahler setting i.e. when $M$ admits an integrable almost
 complex structure $J$ such that the bilinear form
 $g(X,Y)=\omega(X,JY)$ defines a Riemaniann metric on $M$, the
 associated Riemannian properties
of Lagrangian submanifolds have been studied by different authors
(see \cite{HL}, \cite{Ta}, \cite{Br}, \cite{Oh}, \cite{Oh2}),
in particular in relation to the
analysis of {\em minimal} Lagrangian submanifolds. In \cite{Oh}
the author asks for a group theoretical machinery producing
minimal Lagrangian submanifolds in Hermitian symmetric spaces.\\
In the present paper we first study the existence problem of
homogeneous Lagrangian submanifolds in compact K\"ahler manifolds, coming to
the characterization of isometric actions admitting a Lagrangian
orbit, by imposing an additional hypotesis on $M$,
holding for a large class of K\"ahler manifolds
including irreducible Hermitian symmetric spaces.
Namely we require the space $H^{1,1}(M)$ to be
1-dimensional.
\begin{theorem}\label{teoremone} Let $K$ be a compact connected group of isometries
acting in a Hamiltonian fashion on a compact K\"ahler
manifold $M$ with $h^{1,1}(M)=1$. Then $M$ admits a $K$-homogeneous
Lagrangian submanifold if and only if $K^\C$ has an open Stein orbit
in $M$. \end{theorem}
In \cite{Akbook} it is proved that if $G:=K^\C$ acts holomorphically
on a {\it complex} manifold with an open Stein orbit, then there
exists a {\it totally real} $K$-orbit $\mathcal O$, i.e. at every point of
$\mathcal O$ the tangent space does not contain complex lines.\\
A counterexample shows that when $h^{1,1}>1$, even the presence of
an open Stein $K^\C$-orbit does not guarantee the existence of a
Lagrangian $K$-orbit.\\
When $K$ is semisimple it turns out that the Lagrangian orbit is
unique and,  when $M$ is K\"ahler Einstein,  it is also minimal.
In the general case we describe the actions
having infinitely many Lagrangian orbits, characterizing the
minimal ones in Section 2.
\begin{theorem}\label{cardinalita}
A Lagrangian $K$-orbit, $K\cdot p$, is isolated (actually unique)
if and only if the smallest subgroup $K'$ of $K$ such that $K\cdot
p=K' \cdot p$ is semisimple.
\end{theorem}
Our main tool will be the moment map, that can be defined whenever we
consider an Hamiltonian group action on $M$.
More precisely,
let $(M,\omega,J)$ be a compact $2n$-dimensional  K\"ahler manifold, acted on in a Hamiltonian
fashion by a compact connected subgroup $K$ of its full isometry  group.
This means that there exists a smooth map
$\mu: M\rightarrow\mathfrak{k}^*=\rm{Lie}(K)^*$, called a {\em moment map}, with the following properties:
\begin{enumerate}
\item $d\mu_p(v)(X)=\omega_p(v,\widehat {X}_p)$ for all $p\in M$, $v\in T_p M$
and $X\in \mathfrak{k}$. Here $\widehat {X}_p$ stands for the
fundamental field associated to $X$, evaluated at $p$;
\item $\mu$ is $K$-equivariant with respect to the coadjoint action of $K$ on $\mathfrak{k}^*$.
\end{enumerate}
In general the matter of existence and uniqueness of the moment map is delicate.
However, whenever the Lie group $K$ is
semisimple  there is  a unique moment map (see e.g. \cite{Kir}).
If $(M,\omega)$, as in our situation, is a compact K\"ahler
manifold and $K$ is a connected compact group of holomorphic isometries
then the existence problem can be easily solved:
a moment map exists if and only if $K$ acts trivially on the Albanese manifold of $M$ (see e.g. \cite{HW}).
Moreover if $\mu_1$ and $\mu_2$ are two moment maps, there exists $c$ in the dual of the Lie
algebra of the center of $K$, such that
$\mu_1=\mu_2+c$.\\
In \cite{GP} the authors have studied the critical set of the squared moment map
$||\mu||^2$, where $||\cdot||$ denotes the norm induced by an
$Ad(K)$-invariant inner product $<,>$ on $\mathfrak{k}^*$. In particular it is proved that
if a point $x\in M$
realizes the maximum of $||\mu||^2$, then the orbit $K\cdot x$ is
complex; hence $K\cdot x=K^\C \cdot x$ is a closed $K^\C$-orbit; it
is therefore natural to consider the ``dual'' problem, i.e. to investigate
the $K$-orbits through points $y\in M$ that attain
the minimum of $||\mu||^2$.
At least when $K$ is semisimple and $K^\C$ has an open Stein orbit on
$M$, Theorem \ref{teoremone} is a step in this direction.\\
While in Theorem \ref{teoremone} we prove the existence of a Lagrangian orbit $L$,
we do not exhibit an effective way to single out $L$. At least for self-dual representations,
we give in Remark \ref{selfdual} an explicit expression, in terms of the
highest weight vector, of a point through which the orbit
is Lagrangian.
Using this result and several {\it ad hoc} arguments we
finally give the complete classification of Lagrangian
submanifolds of the complex projective space on which a simple
group of isometries of the whole space acts transitively.
\begin{theorem}\label{classificazione}
Let $K$ be a simple compact Lie group acting on the complex projective space $\P(V)$,
by means of a unitary representation
$\rho:K\to \U(V)$.
The group $K$ has a Lagrangian orbit in $\P(V)$  if and only if it appears in Table $1$.
\end{theorem}
\par The paper is organized as follows:
In the first section we introduce some notations and give the proof
of Theorem \ref{teoremone} and Theorem \ref{cardinalita}.
In the second section we analyse
the minimality of Lagrangian submanifolds, while in the last section
we give the complete classification of simple Lie groups that admit an homogeneous
Lagrangian submanifold $L$ in the complex
projective space.\\
{\em Notations and conventions}. Lie groups and their Lie algebras
will be indicated with capital and gothic letters respectively.
Moreover, after identifying, by means of a $Ad(K)$-invariant inner
product $<,>$ on $\mathfrak k ^*,$ the Lie algebra $\mathfrak k$
and its dual ${\mathfrak k}^*$, we will alternatively consider
$\mu$ as a $\mathfrak{k}$-valued map.\\ \noindent {\bf
Acknowledgements.} The authors would like to thank Prof. F.
Podest\`a  for his constant support and his useful advices.
%
%
%
%
\section{Existence and Uniqueness}\label{mainsection}

Let $M$ be a compact complex manifold with a K\"ahler form $\omega$ and $K$ be a compact group of
isometries acting on $M$ in a Hamiltonian fashion.
From now on we fix a moment map $\mu$,
and focus on the set of
points of $M$ sent by $\mu$ to $\mathfrak{z(k)}$, we will
denote this set by $\mathcal Z$.
The defining
properties of $\mu$ imply that the $K$-orbits through points of
$\mathcal Z$ are $\omega$-isotropic, indeed, for every $X,Y$ in
$\mathfrak{k}$ and $q=kp$ in $K\cdot p$
$$\omega_q(\widehat{X}_q,\widehat{Y}_q)=d\mu_q(\widehat{X}_q)(Y)={\frac{d}{dt}_|}_{t=0} \exp t X k\cdot\mu(p)(Y)=0.$$
Since the $K$-action on $M$ is holomorphic, it induces, when $M$ is compact, an action of the complexified group $G:=K^\C$ on $M.$
With these notation we state
\begin{lemma}\label{lemma1} Let $p$ be in $\mathcal Z$.
Then the following statements are equivalent
\begin{enumerate}
\item the $K$-orbit through $p$ is Lagrangian;
\item the $G$-orbit $\Omega$ through $p$ is open in $M$, i.e.
$M$ is a $G$-almost homogeneous space.
\end{enumerate}
In this case the $G$-orbit is a Stein manifold.
\end{lemma}
\begin{proof} Denote by $\mathcal O$ the $K$-orbit through $p$.\\
$(i)\Rightarrow (ii)$ The tangent space to the $G$-orbit through
$p$ is given by
$$T_p (G\cdot p)=T_p \mathcal O+JT_p \mathcal O$$ and
the sum must be direct since $JT_p \mathcal O\cap T_p \mathcal
O=\{0\}$ because $\mathcal O$ is Lagrangian. Hence $\dim T_pG\cdot
p=2\dim T_p \mathcal O =\dim M$ and the $G$-orbit is open.\\
$(ii)\Rightarrow (i)$ Since the $G$-orbit is open, $\dim T_p
\mathcal O+ \dim JT_p \mathcal O\geq 2n.$ Now the conclusion
follows recalling that  $\mathcal O$ is isotropic $(\mu(p)\in
\mathfrak{z(k)})$, hence $\dim \mathcal O \leq n$.\\
 Let $H\leq
K$ be the isotropy subgroup at $p$. When we consider the
complexified action, the Lie algebra of the stabilizer,
$\mathfrak{g}_p$, is given by the set of vectors $W=X+iY$ such
that $\widehat W_p =0$. Now, recalling that $JT_p \mathcal O={(T_p
\mathcal O )}^\perp$ we get that $\widehat {X}_p=\widehat {Y}_p=0$, therefore
the complex isotropy is reductive and the open orbit
$\Omega=\frac{K^\C}{H^\C}$ is Stein thanks to a theorem of
Matsushima \cite{Ma}.
\end{proof}
%
%
%
%
An immediate consequence is the following
\begin{corollary}\label{Cod1}
The complement of $\Omega$ in $M$ has complex codimension $1$.
\end{corollary}
Theorem \ref{teoremone} proves that, by imposing an
additional hypothesis on the cohomology of $M$, the existence
of an open Stein $G$-orbit is indeed sufficient to guarantee the
presence of a Lagrangian $K$-orbit, while in Theorem
\ref{cardinalita} we characterize  the actions having infinitely many Lagrangian orbits.
Now we recall two results that will be used in proving the theorems, the first one
is due to Kirwan \cite{Kir}
\begin{lemma}[Kirwan]\label{Kir} Let $x$ and $y$ be two points in a K\"ahler manifold $M$,
acted on in a Hamiltonian fashion
by a group of isometries $K,$ such that $\mu(x)=\mu(y)=0$. Suppose
that $x$ and $y$ lie in different $K$-orbits, then there exist two $K^\C$-invariant
disjoint neighborhoods $U_x$ and $U_y$ of $x$ and $y$
respectively.
\end{lemma}
The following is a classical result in K\"ahler geometry (see e.g.
\cite{KoMo} for a proof), it is essentially a consequence of
$\partial\overline\partial$-lemma, holding for compact K\"ahler manifolds.
\begin{proposition}\label{kodaira} Let $L\rightarrow M$
be a line bundle on a compact K\"ahler manifold $M$. If $\omega$
is any real, closed $(1,1)$-form such that $[\omega]=c^{\R}_1(L)\in
H_{dR}^2(M)$, then there exists a Hermitian metric along the
fibers of $L$ whose curvature form is $\Theta=\frac{i}{2\pi}\omega$.
\end{proposition}
Now we can  prove Theorem \ref{teoremone}.\\
%
%
%
%
\noindent Note that a compact K\"ahler
manifold $M$ with $h^{1,1}(M)=1$ is necessarily projective. Indeed, since
the K\"ahler form $\omega$ is of type $(1,1)$, we
can scale it  so that we obtain an integral
class $[\widetilde{\omega}]\in H^2(M;\Z)$ and use the Kodaira
Embedding Theorem.
\\
Note that the hypotheses of Theorem \ref{teoremone} are naturally satisfied
when $M$ is a compact irreducible Hermitian symmetric space.

\begin{proof}[Proof of Theorem \ref{teoremone}]
We need only to prove that if $G=K^{\C}$ has an open Stein orbit,
then there exists a Lagrangian $K$-orbit.
Denote again by $\Omega=G\cdot p$ the open Stein orbit and by $Y$
its complement in $M$. By Corollary \ref{Cod1}, $Y$ is a divisor of
$M$ and therefore it determines a holomorphic line bundle $L$ on
$M$ and a section $\sigma\in H^0(L)$ such that $Y$ is the vanishing locus of $\sigma$.
We scale $\omega$ so that we obtain a positive generator of the free part
of $H^2(M;\Z) \cap H^{1,1}(M)$.  Since $h^{1,1}(M)=1$ the first Chern class of $L$ is a positive
integer multiple of the class of the scaled K\"ahler
form $\omega$ on $M$:
\[
c_1(L)=m[\omega]\in H^{1,1}(M).
\]
Now, by Proposition \ref{kodaira}, it is possible to find a Hermitian metric $h$ on the fibers
of $L$ such that its curvature form is
$$\Theta=m\frac{i}{2\pi}\omega.$$
On the other hand, the curvature on $\Omega$ is exactly (see e.g. \cite{KoMo})
$$\partial\overline\partial \log\|\sigma\|^2,$$
where $\|\cdot\|$ is the norm induced by $h$.
Thus we have found a strictly plurisubharminic real valued function
$\rho$
such that $\omega=i\partial\overline\partial \rho.$
Note that
by construction $\rho:\Omega\rightarrow \R$ is an
exhaustion function. Observe that we can assume $\rho$ to be $K$-invariant, since this can be
achieved by averaging over the {\em compact} group $K$.\\
Starting from $\rho$ we can define a map
$\phi:\Omega\rightarrow\mathfrak{k}^*$ as follows
$$\phi(p)(X):=\frac{1}{2}{(J\widehat X)}_p(\rho).$$
Clearly the $K$-invariance of $\rho$ implies the $K$-equivariance of $\phi$.
Moreover, for every $p\in \Omega$, $v\in T_p\Omega$ and $X\in \mathfrak{k}$,
we have that $d\phi_p(v)(X)=\omega_p(v,\widehat X_p)$ (see \cite{HHL} for the proof).
Hence $\phi$ is a moment map for the Hamiltonian action of $K$ on $\Omega$
and therefore its extension to the whole $M$ $\phi$ differs from $\mu$ by an element $z$ of
$\mathfrak{z(k)}$.\\ Let now $x_o\in \Omega$ be a critical point of
the exhaustion function $\rho$, then $d\rho_{x_o}=0$ and
$\phi(x_o)(X)=\frac{1}{2}(J\widehat X_{x_o})(\rho)=0$ for all $X\in \mathfrak{k}$.
Thus $\mu(x_o)$ belongs to the Lie algebra of the center
of $K$, and the $K$-orbit through $x_o$ is Lagrangian by Lemma
\ref{lemma1}.\end{proof}
%
%
%
%
\begin{remark} If the assumption of the Hodge number $h^{1,1}(M)$ in Theorem \ref{teoremone}
is not satisfied we cannot reach the same conclusion. Indeed
consider the example of $\SU(3)$ acting on $\P^2
\times {\P^2}$
as follows
$$A\cdot([x],[y])=([Ax],[\bar{A}y])$$
with $A \in \SU(3)$ and $x,y \in \C^3\setminus \{0\}$.
Since $h^{1,1}(\P^2\times\P^2)=2$ we can choose an $\SU(3)$-invariant
symplectic form $\omega_\varepsilon=\omega_0\oplus (1+\varepsilon)
\omega_0$ on $\P^2\times\P^2$, where $\omega_0$ is the
Fubini-Study $2$-form on $\P^2$ and $\varepsilon$ is a small positive constant.
In this case, there
exists an open Stein $G$-orbit (see e.g. \cite{Ak2}),
while the image of the moment map does not contain $0$ (see also \cite{BG}
for the picture of the moment polytope in this case).
\end{remark}
%
%
%
%
\begin{remark}
If the group $K^\C$ has an open Stein orbit in $M$ with $h^{1,1}(M)=1$,
the same is true for $(K\cdot Z)^\C$, where $Z$ centralizes $K$.
Indeed consider $p \in M$ such that $K\cdot p$
is Lagrangian (cfr. Theorem \ref{teoremone}), then $\mu(p) \in
\mathfrak{z(k)}$ where $\mu$ is a moment map for the $K$-action.
On the other hand $\mu$ is the composition of the moment map $\mu'$ for the action
of $K':=K\cdot Z$ with the projection induced by the inclusion on the dual of the Lie algebras.
Therefore $\mu'(p) \in \mathfrak{z(k')}$ and for dimensional reasons $K'\cdot p$ is Lagrangian
and the claim follows from lemma \ref{lemma1}.
\end{remark}
Under the same assumptions of Theorem \ref{teoremone} on $M$ we
 prove Theorem \ref{cardinalita}
%
%
%
%
\begin{proof}[Proof of Theorem \ref{cardinalita}]
%
%
%
%
In the semisimple case the moment map is unique, therefore, using
the same notation as in  Theorem \ref{teoremone} we have $\mu(x_o)=0$ at
the critical point $x_o$ of $\rho$, and $\mu^{-1}(0)\cap
\Omega\neq\emptyset$. Take $x$ and $y$ in
$\mu^{-1}(0)\cap \Omega$, applying Lemma \ref{Kir}, we deduce
that $x$ and $y$ belong to the same $K$-orbit,
and  $\mu^{-1}(0)\cap \Omega$ is therefore compact. Then, since
when $M$ is  compact the fibers of the moment map are connected \cite{Kir},
we get that $\mu^{-1}(0)$ is contained in $\Omega$ and it is a single $K$-orbit.\\
If the semisimple part of $K$, that will be  denoted by $K_s$, has
a Lagrangian orbit $L$, then $K$ has a unique Lagrangian
orbit. Indeed, combining Theorem \ref{teoremone} and the previous
remark, we get that there exists a Lagrangian $K$-orbit, this is
contained in $\mu_s^{-1}(0)=L$, where $\mu_s$ is the
moment map for the $K_s$ action, and it is therefore unique.\\
%
%
%
%
%
%
%
%
Now assume that $K\cdot p$ is a Lagrangian orbit and denote by $H$
the connected component of the identity of the isotropy subgroup $K_p$. At the Lie algebra level
$\mathfrak{k}$ can be written as the direct sum $\mathfrak{k}_s\oplus\mathfrak{z(k)}$.
Consider the projection $\pi:\mathfrak{k}\rightarrow
\mathfrak{z(k)}$. Suppose that $\pi(\mathfrak{h})\neq 0$, and
consider $Z' \subset Z$ a subtorus such that its Lie algebra satisfies
\[
\mathfrak{z(k)}=\pi(\mathfrak{h})+\mathfrak{z'},
\]
and call $K'$ the group $K_s\cdot Z'$. We first prove that $K'\cdot p$ has the same
dimension of $K\cdot p$ and therefore $K'\cdot p$ is Lagrangian.
The set of tangent vectors to the $K'$-orbit is given by
\[
\widehat {\mathfrak{k}'}_{|p}=\widehat {\mathfrak{k}_s}_{|p}+\widehat
{\mathfrak{z}'}_{|p}
\]
while the set of vectors tangent to the
$K$-orbit is given by
\[
\widehat {\mathfrak{k}}_{|p}=\widehat {\mathfrak{k}_s}_{|p}+\widehat{
\mathfrak{z}'}_{|p}+\widehat {\pi(\mathfrak{h})}_{|p}.
\]
By construction
$\widehat {\pi(\mathfrak{h})}_{|p}$ is contained in $\widehat
{\mathfrak{k}_s}_{|p}$ hence the $K-$ and $K'$-orbits through $p$ coincide.\\
Now denote by $H'$ the group ${K'_p}^o$; the
projection $\pi(\mathfrak{h'})$ is $\{0\}$, i.e. $H'\subseteq K_s$. \\
%
%
%
%
We claim that for all subtori $Z''\leq Z'$ of codimension $1$ in $Z'$, the
group $K''=K_s\cdot Z''$ has no Lagrangian orbits. This can be
proven observing that $H'$ coincides with ${K''_p}^o$.
Indeed $K''\cdot p\subset K'\cdot p$ and $\rm{codim}_{K'\cdot
p}(K''\cdot p)=1$, hence $$\dim K'_p=\dim K'-\dim K'\cdot p=\dim
K''-\dim K'\cdot p-1=\dim K''-\dim K''\cdot p=\dim K''_p$$
therefore
$$\dim K''\cdot p< \dim K'\cdot p$$
and the $K''$-orbit is not Lagrangian.\\
%
%
%
%
Denote by $\mu''$ the moment map associated to the $K''$-action on
$M$. Consider the set $M^{{K'_p}^o}$, i.e the set $$\{x\in
M|\;{H'}\cdot x=x\}.$$ We first prove that
$\mu''(M^{{H'}})$ is contained in
$\mathfrak{z}_{\mathfrak{m''}}(\mathfrak{h'})$, where
$\mathfrak{k}''=\mathfrak{h'}\oplus \mathfrak{m''}$ and
$\mathfrak{z}_{\mathfrak{m''}}(\mathfrak{h'})=\{X\in
\mathfrak{m''}|\;[X,\mathfrak{h'}]=0\}.$
Clearly
$\mu''(M^{{H'}})$ is contained in
$\mathfrak{z}_{\mathfrak k''} (\mathfrak{h'})$. Moreover let $\gamma(t)$ be
a smooth curve  contained in $M^{{H'}}$ joining $p$ and a
point $x\in M^{{H'}}$. We get
$$
\frac{d}{dt}<\mu''(\gamma(t)),\mathfrak{h'}>=<d\mu''_{\gamma(t)}(\gamma'(t)),\mathfrak{h'}>
=\omega_{\gamma(t)}(\gamma'(t),\widehat{\mathfrak{h'}}{_{|_{\gamma(t)}}})\equiv0
$$
 where the last equality holds since
$\widehat{\mathfrak{h'}}_{|\gamma(t)}= 0.$
Now recall
that the orbit $K'\cdot p$ is Lagrangian hence $\mu''(p)=c\in
\mathfrak{z(k'')}=\mathfrak{z''}$ which is orthogonal to
$\mathfrak{h'}$. Therefore $\mu''(x)$ is orthogonal to
$\mathfrak{h'}$ and belongs to $\mathfrak{m''}\cap
\mathfrak{z}_{\mathfrak{k''}}(\mathfrak{h'})$=$\mathfrak{z}_{\mathfrak{m''}}(\mathfrak{h'})$
as claimed.\\
%
%
%
%
Now the dimension of $M^{{H'}}$ is given by
$$
\dim M^{H'}=2\dim {(K'\cdot p)}^{H'}=2(\dim {(K''\cdot
p)}^{H'}+1)= 2(\dim\mathfrak{z}_{\mathfrak{m''}}(\mathfrak{h'})+1).
$$
This will be used in proving that $Q=\mu''^{-1}(c)\cap M^{{H'}}$ is a
submanifold.\\Note that
$$
\Ker\; d\mu''_p=\{Y\in T_p M|\;
\omega(Y,\widehat{X}_p)=0 \; \mbox{for all} \; X\in \mathfrak{k''}\}={(
J\widehat {\mathfrak{k''}}_{|p})}^\perp=T_p K'\cdot p\oplus V_1
$$
where $V_1$ has dimension $1$ and is contained in ${T_p K'\cdot
p}^\perp$, indeed $K'\cdot p$ is Lagrangian and $T_p M=T_p K'\cdot
p\oplus  J\widehat {\mathfrak{k''}}_{|_p}\oplus V_1$. Moreover $V_1$ is
contained in $T_p M^{{H'}}$; indeed ${H'}$ acts by
isotropy on $T_p M$ and leaves $\widehat {\mathfrak{k''}}_{|_p}$ invariant
hence ${( J\widehat {\mathfrak{k''}}_{|_p})}^\perp$ invariant. Therefore
$V_1$ is ${H'}$-invariant, hence fixed, since it is $1$
dimensional and ${H'}$ is compact.\\
We conclude that
$$
\Ker \;d\mu''_p\cap T_p (M^{H'})= {(T_p
K'\cdot p)}^{{H'}} \oplus V_1
$$
and
$$
\dim ({\Ker \; d\mu''_p\cap T_p (M^{H'})})=\dim \mathfrak{z}_{\mathfrak{m''}}(\mathfrak{h'})+2.
$$
Counting the dimension of the image, it follows that
${\mu''}_{|M^{H'}}$ is a submersion at $p$, hence $Q$ is a
manifold locally around $p$ whose dimension is $\dim
\mathfrak{z}_{\mathfrak{m''}}(\mathfrak{h'})+2$. Note that for dimensional reasons $Q\setminus
K'\cdot p\neq \emptyset$. To complete the proof it
is sufficient to observe that if we take $y\in Q\setminus K'\cdot
p$, sufficiently close to $p$, then the $K'$-orbit through $y$ is
Lagrangian. Indeed $K'\cdot y$ is isotropic for $\mu''(y)\in \mathfrak{z''}$, and furthermore
${H'} \subseteq K'_y$; by the Slice Theorem ${K'_y}^o$ is conjugated to a subgroup of $H'$ hence
$\dim K'_y \leq \dim H'$ so that $\dim K'\cdot y=\dim K' \cdot p$.
\end{proof}

The uniqueness of the Lagrangian $K$-orbit in the semisimple case is
independent of the assumption $h^{1,1}(M)$.\\
Note that whenever  $\mu^{-1}(0)\cap \Omega\neq\emptyset$,
we can argue that $\Omega$ coincides the
set of semistable points $M^{ss}:=\{x\in M| \: \overline{G\cdot x}\cap
\mu^{-1}(0)\neq\emptyset\}$.
Indeed $\Omega$ is always contained in $M^{ss}$, moreover (see e.g. \cite{Sj})
$M^{ss}$ is the smallest $G$-invariant subset of $M$ that contains $\mu^{-1}(0)$,
therefore $M^{ss}=\Omega$. In the K\"ahler case it is easy to see that
the stratum associated to the minimum critical set
of $\|\mu\|^2$ (see \cite{Kir} for the precise definition) coincides with the set
of semistable points.
Lerman in \cite{Le} shows that this stratum retracts to the zero set; we have thus proved
the following
\begin{proposition} If  $\mu^{-1}(0)\cap \Omega$ is not empty, then
$\Omega$ has  $\mu^{-1}(0)$ as a
deformation retract, and thus has the same homotopy type of the Lagrangian orbit.
\label{retraz}
\end{proposition}
\begin{remark}
Let $K\cdot p$ be a Lagrangian $K$-orbit. One can easily show
that, if $Z\in \mathfrak{z(k)}$, then $K_p=K_{\exp iZp}$ and all
the orbits through $\exp iZ p$ are totally real ({\it i.e.} at every point the
tangent space is transversal to its image via the complex structure), but not in general
Lagrangian. Consider for example the action on $\C\P^N$,
with $N={\frac{n^2+3n}{2}-1}$ induced by the representation $\rho$  of $T^2\times SU(n)$
on $V=S^2(\C^n)\oplus \C^n$,
defined by $\rho(g)(X,Y)=(\alpha AXA^t,\beta A y)$
for $g=(\alpha,\beta,A)\in T^2\times \SU(n)$,
where we see the elements of $S^2\C^n$ as symmetric matrices. Here there are more than one
$K$-Lagrangian orbit but,
 moving through points $\exp i tZ p$, one does not meet any (other)
 Lagrangian orbit.
Azad, Loeb and Qureshi in \cite{ALQ} give necessary and sufficient
conditions under which one can prove that there are infinitely many
totally real orbits; more precisely this is the case whenever
$N_G(G_p)/G_p$ is not finite. In the non semisimple case this
condition is always satisfied.
\end{remark}
\begin{remark} Whenever the isotropy of a Lagrangian $K$-orbit is
discrete, the set of Lagrangian orbits is a manifold whose
dimension equals the dimension of the center of the group
\cite{Pa}. This situation holds whenever there exists a regular
(i.e. principal or exceptional) Lagrangian $K$-orbit. Nevertheless
note that if $p$ belongs to the set
\[
M_{\mu}\cap M_{princ}
\]
where $M_{\mu}$ is
the set of points $x$ in $M$ whose orbits $K\cdot \mu(x)$ has
maximal dimension and $M_{princ}$ the set of principal points in
$M$, and the $K$-orbit through $p$ is Lagrangian, then necessarily
$K$ must be abelian. Indeed, in general when $p\in M_{\mu}\cap M_{princ}$,
$K_{\mu(p)}/K_p$ is abelian (see e.g.\cite{HW}); since in this case $K\cdot p$ is principal and
Lagrangian $K_p$ is trivial, $K_{\mu(p)}$
is abelian, but $\mu(p)\in \mathfrak{z(k)}$, hence $K_{\mu(p)}=K$
and the claim follows.
\end{remark}
\begin{remark}
As a consequence of Theorem \ref{teoremone}
we have incidentally proved that linear representations of
complex semisimple Lie groups are {\em balanced}, in the sense of \cite{Wi},
if and only if they have an open Stein orbit.
\end{remark}
\section{Minimality of Lagrangian orbits}
We here give a proof of the minimality of the Lagrangian orbit
in the semisimple case, however this can be proved also as a
consequence of the more general fact stated in Proposition \ref{minimale}.
\begin{proposition} If $K$ is semisimple and $M$ is K\"ahler
Einstein, the $K$-orbit is also minimal.
\end{proposition}
\begin{proof} If $H$ denotes the mean curvature vector of the Lagrangian orbit $\mathcal{O}$, it is known (see
Dazord \cite{Da}) that the $1$-form $\alpha\in\Lambda^1(\mathcal O)$
which is the $\omega$-dual of $H$ restricted to $\mathcal O$ is
closed. But $\alpha$ is $K$-invariant, hence for every $X,Y\in
\mathfrak{k}$ we have
$$0=d\alpha(\widehat X,\widehat Y)=\widehat X\alpha(\widehat Y)-\widehat Y \alpha
(\widehat X)-\alpha([\widehat X,\widehat Y])=-\alpha([\widehat
X,\widehat Y]),$$ so that
$\alpha([\widehat{\mathfrak{k}},\widehat{\mathfrak{k}}])=\alpha(\widehat{\mathfrak{k}})=0$
and $\alpha\equiv 0$. This means $H=0$.
\end{proof}
Actually one can characterize the minimal Lagrangian orbit $L$
in the general case.\\
When $(M,\omega)$ is compact  we can define a {\it canonical moment map}, $\widetilde{\mu}$, that is characterized by
the fact that $\int_M \mu \omega^n=0$. If further $M$ is K\"ahler-Einstein with Einstein constant $c$, then
$\widetilde{\mu}$ can be explicitly written (see e.g. \cite{Fu},
\cite{Po}):
$$\widetilde{\mu}(p)(Y):=\frac{1}{2c}\mbox{div} (J\widehat{Y}_p)$$
for every $Y\in \mathfrak{k}.$
\begin{proposition} \label{minimale} Let $\widetilde\mu$ be the canonical moment 
map of a K\"ahler-Einstein manifold,
then a Lagrangian orbit $\mathcal O$ is minimal if and only if $\widetilde\mu(\mathcal O)=0.$
\end{proposition}
The previous result is stated and proved in \cite{Pa} assuming that the Lagrangian orbit is {\it principal}.
Actually Proposition \ref{minimale} holds without any assumption on the type of Lagrangian orbits.
Indeed, since $L$ is Lagrangian, in order to prove that $L$ is minimal,
it is sufficient to show that the mean curvature vector $H$ at some point $p$ of $L$
is orthogonal to $J\widehat{\mathfrak{k}}_p$ as done in \cite{Pa} in Proposition $5$.
Once an orthonormal frame $\{e_i\}$ at
$p$ is fixed, we have
\begin{eqnarray*}
\langle H, J\widehat{Y} \rangle & = & \langle \nabla_{e_i}e_i,J\widehat{Y}\rangle\\
                                & = & e_i \sum\langle e_i , J\widehat{Y} \rangle -
                                          \sum\langle e_i,\nabla_{e_i}J\widehat{Y} \rangle\\
                                & = & -\sum \langle e_i,\nabla_{e_i}J\widehat{Y} \rangle\\
                                & = & -\frac{1}{2} \mbox{div}J\widehat{Y}\\
                                & = & c\widetilde{\mu}_p(Y) = 0.
\end{eqnarray*}


Combining the previous proposition and the fact that the zero level set of the 
moment map is a single orbit  when it meets the
open Stein $K^\C$-orbit (see proof of Theorem \ref{cardinalita}), we get

\begin{corollary}\label{unica}
Let $K$ be a compact connected group of isometries
acting in a Hamiltonian fashion on a compact K\"ahler-Einstein
manifold $M$. Then $M$ admits at most one minimal Lagrangian $K$-orbit. 
\end{corollary}
\begin{corollary} Under the same hypotheses of Theorem \ref{teoremone}, assuming further 
that $K^\C$ is simply connected, $M$ is K\"ahler-Einstein and the 
isotropy subgroup at a point of the Stein orbit has finite connected 
components we get that $M$ admits  a unique $K$-orbit wich turns out to be minimal.
\end{corollary}
\begin{proof} From Theorem \ref{teoremone} we get that there is a 
Lagrangian $K$-orbit $L$; moreover any other Lagrangian $K$-orbit 
has the same homotopy type of $L$ by Proposition \ref{retraz} and 
therefore has finite fundamental group. But, according to Chen 
(see Theorem 5.1 in \cite{Chen} and the reference therein), 
in a K\"ahler-Einstein manifold, the mean curvature of
 every compact Lagrangian submanifolds with $b_1=0$ must vanish somewhere. 
The homogeneity implies 
that all the Lagrangian orbits are minimal. 
The conclusion follows from Corollary \ref{unica}.
\end{proof}
 Obviously the same result holds if $K^\C$ is only supposed to have 
finite fundamental group.
 
\section{The classification of Simple Lie groups with a Lagrangian orbit 
in the complex projective space}
In this section we give the complete  classification
of simple compact  Lie groups $K$ with a Lagrangian
orbit in the complex projective space.
We give also an explicit description of Lagrangian orbits,
except in case $K=E_7$.
This part can be treated combining the results of section
\ref{mainsection} with the work of Sato and Kimura \cite{SK} and Kimura \cite{Ki}.\\
Consider a finite-dimensional unitary representation of a compact Lie
group $K$ on a Hermitian vector space $(V,\langle,\rangle)$. Endow
$\P(V)$ with the Fubini-Study K\"ahler form and consider the induced
$K$-action. Note that this action is automatically Hamiltonian since
$\P(V)$ is simply connected. The map $\mu: \P(V) \to \mathfrak{k}^*$
defined for every $v \in V$  and $X \in \mathfrak{k}$ by
\begin{equation}
\label{mmproj}
\mu([v])(X)=\frac{1}{i}\frac{\langle X\cdot v,v\rangle}{\langle v, v \rangle}
\end{equation}
is a moment map for the $K$-action on $\P(V)$.  \\
Here we recall notations and results from \cite{SK}. Given  a connected complex linear algebraic group $G$,
and  a rational representation $\rho$ of $G$ on a finite dimensional complex vector space $V$,
a triplet $(G,\rho,V)$ is {\em prehomogeneous} if $V$ has a Zariski dense $G$-orbit.\\
We give here an easy-to-prove lemma that allows to find relations  between almost homogeneous
actions on the projective space and  prehomogeneus triplets.
\begin{lemma} \label{sato}
Let $G$ be any complex, connected Lie group. $G$ acts with an open
dense orbit on $\C\P ^{n-1}$ if and only if $G \times GL(1)$ acts
with an open dense orbit on $\C^n$ i.e. $(G\times
\GL(1),\rho,\mathbb{C}^n)$ is a prehomogeneus triplet.
\end{lemma}
Hence, thanks to Theorem \ref{teoremone}, in order to classify
the action of compact {\em simple} Lie groups on the projective space
admitting a Lagrangian orbit,
it is sufficient to go through the list of prehomogeneous triplets in \cite{SK}, and
consider those that have reductive generic isotropy, i.e. those
that have an open Stein $G$-orbit. They are exactly {\it regular} PV
spaces of \cite{SK} (p. 59). These spaces are characterized by the existence
of a {\em relative invariant}, i.e. a rational function $f$ such that
there exists a rational character $\chi$ of $G$ satisfying
$f(\rho(g)x)=\chi(g)f(x)$ for any $g\in G$ and $x\in V$.
We here enclose a lemma that will be useful in the sequel;
the proof can be found in \cite{SK} p.64 .
\begin{lemma} If $\rho$ is an irreducible representation,
 then the polynomial $f$ that defines the hypersurface $Y$ is irreducible.
\end{lemma}

In \cite{SK} prehomogeneous vector spaces are classified up to an equivalence
relation which we are going to describe.

\begin{definition} Two triplets $(G,\rho,V)$ and $(G',\rho',V')$ are called
{\em equivalent} if there exist a rational isomorphism
$\sigma:\rho(G)\to\rho'(G')$ and an isomorphism $\tau:V\to V'$,
both defined over $\C$ such that the diagram
\[
\xymatrix{ V\ar[r]^-\tau \ar[d]_-{\rho(g)} & V'
\ar[d]^-{\sigma(\rho(g))} \\ V\ar[r]^-\tau    &   V'}
\]
is commutative for all $g\in G$. This equivalence relation will be
denoted by $(G,\rho,V)\cong (G',\rho',V')$.
\end{definition}
We say that two triplets $(G,\rho,V)$ and $(G',\rho',V')$ are {\em
castling transforms} of each other when there exist a triplet
$(\tilde{G},\tilde{\rho},V(m))$ and a positive  number $n$ with
$m>n\geq 1$ such that $$(G,\rho,V)\cong (\tilde{G}\times \SL(n)
,\tilde{\rho}\otimes \Lambda_1,V(m)\otimes V(n))$$ and
$$(G',\rho',V')\cong (\tilde{G}\times \SL(m-n)
,{\tilde{\rho}}^*\otimes \Lambda_1,V(m)^*\otimes V(m-n)),$$ where
${\tilde{\rho}}^*$ is the dual representation of
${\tilde{\rho}}$ on the dual vector space $V(m)^*$ of $V(m)$. We
recall that $V(n)$ is a complex vector space of dimension $n.$ A
triplet $(G,\rho,V)$ is called {\em reduced} if there is no
castling transform $(G',\rho',V')$ with $\dim V'<\dim V.$\\
Note that in fact in each class there is only one representative of
the form $G \times \GL(1)$ where $G$ is simple and it is necessarily reduced.
This can be seen {\it a posteriori} as follows.
Suppose that $(G',\rho',V')$ is a reduced and castling
equivalent to  $(G \times \GL(1),\rho,V)$, then there should exist a
representation $\widetilde{\rho}: \widetilde{G} \to \GL(V(m))$ such that
$(G\times\GL(1),\rho,V)\cong (\tilde{G}\times \SL(n)
,\tilde{\rho}\otimes \Lambda_1,V(m)\otimes V(n))$. But now we would
have (at least locally) $G\simeq\SL(n)$ and $\GL(1)=\widetilde{G}$ since
$G$ is simple, hence $G'=\GL(1)\times\SL(m-n)$, but the correspondent
triple does not appear in the list of \cite{SK} (p. 144--146).

\subsection{Stabilizer and fundamental group}
We here collect some results and remarks that will be used in order to single out Lagrangian
homogeneous submanifolds in the complex projective space.\\
Assume that a complex Lie group $G=K^\C$ acts with an open Stein
orbit $\Omega=\P^n\setminus Y$ on $\P^n$.
Denote by $L$ the Lagrangian  $K$-orbit.
Thanks to Proposition \ref{retraz}, we get that $\Omega$ has the same homotopy type of $L$.
We give here a well known result on the topology of the complement of an algebraic
hypersurface $Y$ in $\P^n$ (see e.g.\cite{Li}):
\begin{proposition}
Let $Y$ be an algebraic hypersurface of $\P^n$. If its
irreducicible components $Y_1,\ldots,Y_r$ have degree $d_1,\ldots,d_r$
respectively, then $H_1(\P^n\setminus Y;\Z)=\Z^r/(d_1,\ldots,d_r)$.
\end{proposition}
From the previous proposition it follows that, if $Y$ is irreducible
of degree $d>1$, then $H_1(\P^n\setminus Y;\Z)$ is cyclic of order $d$. The open Stein orbit 
$\P^n\setminus Y$ contains the Lagrangian orbit $L=K/K_p$ and retracts
onto it. From the homotopy sequence, whenever $K$ is simply connected
$\pi_1(K/K_p)\simeq K_p/K_p^o.$ Hence we get a method in order to determine the
number of connected components of the stabilizer $K_p$.
If $N_K(K_p^o)/K_p^o$ is abelian then
$K_p/K_p^o=\Z_d$.
Indeed $K_p/K_p^o\subset N_K(K_p^o)/K_p^o$ is abelian, hence
$$K_p/K_p^o=\pi_1(K/K_p)=H_1(K/K_p)=\Z_d.$$
\begin{remark} ({\it Self-dual representations})\label{selfdual} Let
$V$ be a $(N+1)$-dimensional complex self-dual representation of a
compact Lie group $K$ and $\mu$ be the corresponding moment map.
Assume that $G=K^\C$ has an open Stein orbit $\Omega=G/H$ in $\P^N$.
Assume also that the highest weight $\lambda$ of the representation satisfies
$2\lambda \notin R^+$. Denote by ${\mathcal{P}}=-{\mathcal{P}}$ the
set of weights.
If $v_{\pm 1}\in V_{\pm \lambda}$ are two non zero vectors with the
same norm, then $[v]:=[v_1+v_{-1}] \in \P^N$ is a point in
$\mu^{-1}(0)$ (see \cite{DK}). If moreover $2\lambda \notin R^++R^+$, then
\[
(\mathfrak{k}_{[v]})^\C=\ker \lambda \oplus_{\pm\alpha\in A_\lambda}\mathfrak{k}_\alpha,
\]
where $A_\lambda=\{\alpha \in R^+: -\lambda +\alpha \notin
{\mathcal{P}}\}=\{\alpha \in R^+: \langle\lambda,\alpha\rangle=0\}$. Indeed
\[
X=H+\sum_{\alpha \in \R^+}c_\alpha E_\alpha+\sum_{\alpha \in \R^-}d_\alpha E_\alpha
\] belongs to $(\mathfrak{k}_{[v]})^\C$ if and only if
\[
X\cdot v =\lambda(H)(v_1-v_{-1})+\sum_{\alpha\in R^+}c_\alpha E_\alpha
v_{-1}+\sum_{\alpha\in R^-}d_\alpha E_\alpha v_1=c\cdot v
\]
and the conclusion follows from the fact that the weight spaces $V_{-\lambda+\alpha}$ and
 $V_{\lambda-\beta}$ are distinct for $\alpha, \beta \in R^+$.
\end{remark}
\begin{remark}If $\Omega=K^\C/H^\C$ is the open Stein orbit, then
  there exists $p\in \Omega$ such that $K_p=H$. Now, by the
  $K$-equivariance of $\mu$, $H=K_p\subseteq K_{\mu(p)}$ which is the
  centralizer of a torus $T$ in $K$.\\ In some situation the only
  centralizer of a torus which contains $H$ is the whole group $K$. In
  this case we have $K_{\mu(p)}=K$ and we can conclude that $\mu(p)=0$
  if $K$ is semisimple.
\end{remark}
\subsection{The case-by-case classification}
In what follows a compact Lie group $K$ acts on the complex finite-dimensional
vector space $V$ by a linear representation $\rho$. Moreover we will identify the
fundamental highest weights $\Lambda_l$ with the corresponding irreducible representations.
\begin{enumerate}
\item \label{2L1} $K=\SU(n),\rho=2\Lambda_1$. Identify the representation space
$V$ with the set of symmetric $n$ by $n$ complex matrices. Now the
Hermitian product on $V$ preserved by $K$ is explicitly given by
$\langle A, B\rangle=\tr (A\overline{B})$ and we get immediately $\mu(I_n)=0$.
Moreover if $Q$ is the $n$ by $n$ matrix $\mbox{diag}(-1,1,\dots,1)$, the stabilizer at $I_n$ is
\[
\{\alpha\cdot \SO(n): \alpha^n=1\} \cup
\{\alpha Q\cdot \SO(n):\alpha^n=-1\}
\]
Therefore the $K$-orbit through $I_n$ is
Lagrangian in $\P(V)$ and $K_{I_n}/K_{I_n}^o\simeq Z_{n}$.
Indeed it is generated by $e^{i\frac{\pi}{n}}$ if $n$ is even, and by
$e^{i\frac{2\pi}{n}}$ if $n$ is odd.
\item {\it $K=\SU(n)$, $\rho=\Lambda_1\oplus \Lambda_1^*$}.
 Identify $V$ with $\C^n\oplus {\C^n}^*$. Take $p=(e_1,e_1^*)$. A direct calculation shows that $\mu(p)=0$.
 The real isotropy is $\SU(n-1)\cdot \Z_2$.
\item {\it $K=\SU(n)$, $\rho=\Lambda_1\oplus\cdots\oplus\Lambda_1$ $n$ times}.
 Identify $V$ with $\C^n\oplus\cdots\oplus\C^n$. Take $p=(e_1,e_2,\ldots,e_n)$. A slightly more
complicated calculation shows that $\mu(p)=0$.
 The complex isotropy of $p$ is discrete while the real one is $\Z_n$.
\item $K=\SU(2n)$, $\rho=\Lambda_2$. Identify the representation space $V$ with the set of
anti-symmetric $2n$ by $2n$ complex matrices. The argument of case
\ref{2L1} applies to
$p=J_n=\left[
\begin{array}{cc}
0 & -I_n \\
I_n & 0
\end{array}
\right]$.
The real stabilizer is
\[
\{\omega\cdot \Sp(n): \omega^{4n}=1\}.
\]
Since $-I_{2n}\in \Sp(n)$ we have $K_{J_n}=\Sp(n)\cdot\Z_{2n}$ and
the $K$-orbit through $J_n$ is Lagrangian.
\item {\it $K=\SU(2n+1)$, $\rho=\Lambda_2\oplus\Lambda_1$}.
Identify the $\Lambda_2$ part of $V$ with anti-symmetric complex
matrices and take $p=(\widetilde{J}_n,e_1)$ where
$\widetilde{J}_n=\left[
\begin{array}{cc}
1 & 0 \\
0 & J_n
\end{array}
\right]$. Again, if $\mu$ is the moment map associated to the
hermitian metric $h((X,v),(Y,w))={\mbox{Tr}}(^t\!X
\overline{Y})+2 ^t\!v\overline{w}$, a straightforward computation
proves that $\mu(p)=0$, and the real isotropy at $p$ is
$\Sp(n)\Z_{n+1}$.
\item
$K=\SU(2)$, $\rho=3\Lambda_1$.
This case has also
been treated in \cite{Ch}. The representation is self-dual, hence
we apply remark 7. Here $\lambda=3\epsilon_1$ and the
set of simple roots $R=\{\pm\alpha\}$ with
$\alpha=\epsilon_1-\epsilon_2.$ Hence $\mathcal P$ is  the set
$\{\lambda,\lambda-\alpha,\lambda-2\alpha,\lambda-3\alpha\}$ and
$\mathfrak{k}_{[v]}=\{0\}$. Explicitly, identifying the representation
space with the space of complex homogeneous polynomial of degree 3,
we can take $[v]=z_1^3+z_2^3$ and $K_{[v]}$ is a non-abelian group of
order 12 whose abelianization is isomorphic to $\Z_4$.
More precisely $K_{[v]}$ is isomorphic to the unique non-trivial
semidirect product $\Z_3\rtimes \Z_4$ in which $\Z_3$ is normal.
\item $K=\SU(6)$ $\rho=\Lambda_3$. The representation is again self-dual, here
$\lambda=\epsilon_1+\epsilon_2+\epsilon_3$ and $\mathcal
P=\{\epsilon_i+\epsilon_j+\epsilon_k; i<j<k\}$ and
$A_\lambda=\{\epsilon_i-\epsilon_j;
i<j<3\}\cup\{\epsilon_i-\epsilon_j;4\leq i<j\}$ hence
$\mathfrak{k}_{[v]}=\mathfrak{su}(3)\oplus\mathfrak{su}(3)$.
Explicitly $[v]=[e_1\wedge e_2\wedge e_3+e_4\wedge e_5\wedge e_6]$
and $K_{[v]}$ has four connected components given by
\[
\left \{ \left[
\begin{array}{cc}
A & 0 \\
0 & D
\end{array}
\right]:
A,D\in \SU(3) \right \} \cup
\left \{ \left[
\begin{array}{cc}
A & 0 \\
0 & D
\end{array}
\right]:
A,D\in \U(3), \det A=\det D=-1\right\}
\]
\[
\cup
\left \{ \left[
\begin{array}{cc}
0 & B \\
C & 0
\end{array}
\right]:
\det B=\det C=i \right \}\cup\left \{ \left[
\begin{array}{cc}
0 & B \\
C & 0
\end{array}
\right]:
\det B=\det C=-i \right \}
\]
Hence the fundamental group of the Lagrangian orbit has
order $4$. But,since $H_1(L,\Z)$ is equal to $\Z_4$ (indeed the invariant has
degree $4$ \cite{SK} (p.144)), $\pi_1(L)=\Z_4$.
\item \label{SU(7)} $K=\SU(7)$, $\rho=\Lambda_3$. Take $p$ such that
$K\cdot p$ is the Lagrangian $K$-orbit in $\P(V)$. By \cite{SK} (p. 144) we
know that $K_p^o=\G_2$. Let $g\in N_K(\G_2)$, then $g$
induces an automorphism of the Lie algebra $\mathfrak{g}_2$ which
is necessarily inner, since $\mathfrak{g}_2$ has only inner
automorphisms. Therefore there exists $h\in \G_2$ such that $gh$
induces the identity on $\mathfrak{g}_2$, i.e. centralizes $\G_2$.
Now, recalling  that $\G_2$ acts irreducibly on $\C^7$, we get
that $gh$ is a scalar multiple of the identity and
$N_K(\G_2)\subset \G_2\cdot \Z_7,$ where $\Z_7$ is the center of
$\SU(7)$, and $K_p=\G_2\cdot \Z_7$.
\item $K=\SU(8)$, $\rho=\Lambda_3$.
In this case, if $p$ is such that
$K\cdot p$ is the Lagrangian $K$-orbit in $\P(V)$, following
the explicit calculations in \cite{SK} (p.87--90), we know that $K_p^o$ is
the image in $\SU(8)$ of $\SU(3)$ via the map
$Ad^{\C}:\SU(3)\rightarrow \Aut(\mathfrak{sl}(3,\C))$, hence
$K_p^o\simeq \SU(3)/\Z_3$. We claim that the cardinality of
$N_K(K_p^o)/K_p^o$ is not greater than $16$, therefore the
cardinality of $H_1(K/K_p,\Z)$ cannot be greater than $16$, while
from \cite{SK} we know that its cardinality is exactly $16$.
Recall that every automorphism of $\mathfrak{su}(3)$ is given by
the composition of an inner and an outer (the conjugation
$\sigma$) automorphism; let $g$ be in $N_K(K_p^o)$ and $\phi_g$
the induced automorphism on $K_p^o$. Then two possibilities arise.
In the first case there exists $h\in K_p^o$ such that
$\phi_g=\phi_h$,
 in other words $gh^{-1}$ commutes with $K_p^o$, which acts irreducibly on $\C^8$,
hence, by the Schur Lemma, it is a scalar multiple of the
identity, i.e. an element of the center $\Z_8$ of $\SU(8)$.
Otherwise  there exists $h\in K_p^o$ such that $\phi_g=\phi_h\circ
\sigma$; in this case put $g_o=h^{-1}\circ g$. Therefore
$N_K(K_p^o)=K_p^o\cdot(\Z_8\cup g_o\Z_8)$, and has at most order
$16$. Now, since $K_p^o$ has no center, we conclude that
$K_p=K_p^0\cdot\Z_{16}$.
\item $K=\Sp(n)$, $\rho=\Lambda_1\oplus\Lambda_1$.
 Identify $V$ with $\C^n\oplus \C^n$. Take $p=(e_1,e_2)$, $\mu(p)=0$. The
 complex isotropy at $p$ is locally isomorphic to $\Sp(n-1,\C)$ while the real isotropy
 is $\Sp(n-1)\cdot \Z_2$.
\item $K=\Sp(3)$, $\rho=\Lambda_3$. The action is
the restriction of the $\SU(6)$ action on the same space.
Therefore the stabilizer is given by the intersection of $\Sp(3)$
with the stabilizer obtained in $(7)$. Hence $K_{[v]}$ is
\[
\left \{ \left[
\begin{array}{cc}
A & 0 \\
0 & \overline{A}
\end{array}
\right]:
A\in \SU(3) \det A=1\right \} \cup
\left \{ \left[
\begin{array}{cc}
A & 0 \\
0 & \overline{A}
\end{array}
\right]:
A\in \U(3), \det A=-1\right\} \cup
\]
\[
\left \{ \left[
\begin{array}{cc}
0 & B \\
-\overline {B} & 0
\end{array}
\right]:
B\in \U(3); \det B=i \right \}\cup\left \{ \left[
\begin{array}{cc}
0 & B \\
-\overline{B} & 0
\end{array}
\right]:
B\in \U(3); \det B=-i \right \}
\]
And we conclude as in $(7)$.
\item {\it $K=\SO(n)$, $\rho=\Lambda_1$}. The representation $\rho$ is self-dual,
nevertheless it is easier to see that $\mu(p)=0$, where $p=[1:0:\dots:0]$.
and $K_p=\SO(n-1)\cdot\Z_2$.
\item {\it $K=\Spin(7)$, $\rho=\mbox{spin rep.}$}
\label{Spin(7)} The orbits of $\Spin(7)$ are the same of $\SO(8)$
(see the previous case), therefore the Lagrangian orbit is
\[
\frac{\Spin(7)}{\G_2\cdot\Z_2}=\frac{\SO(8)}{\SO(7)\cdot\Z_2}=\R\P^7.
\]
\item {\it $K=\Spin(9)$, $\rho=\mbox{spin rep.}$} The case is completely analogous to
the previous one considering the inclusion $\Spin(9) \subset \SO(16)$. Thus the Lagrangian orbit is
\[
\frac{\Spin(9)}{\Spin(7)\cdot\Z_2}=\frac{\SO(16)}{\SO(15)\cdot\Z_2}=\R\P^{15}.
\]
\item $K=\Spin(10)$, $\rho=\Lambda_e\oplus\Lambda_e$ where $\Lambda_e$ is the even
half-spin representation.
The complex isotropy through the point $p=(1+e_{1234},e_{15}+e_{2345})$
is locally isomorphic to $G_2^\C$ (see \cite{SK} also for notations and conventions
on the spin representation space). Moreover a direct computation using formula \ref{mmproj}
shows that $\mu(p)=0$, thus $p$ belongs to a Lagrangian orbit.
\item {\it $K=\Spin(11)$, $\rho=\mbox{spin rep.}$}
This case and the next one (to which we refer) can be treated simultaneously
since $\Spin(11)$ and $\Spin(12)$ have the same orbits on $\P^{31}$.
This can be easily seen noting that
$\Spin(11) \subset \Spin(12)$ and computing the cohomogeneity
of these actions.
In the case of $\Spin(11)$ the isotropy of the Lagrangian orbit
is locally isomorphic to $\SU(5)$.
\item {\it $K=\Spin(12)$, $\rho=\Lambda_e$}
The computation of the fundamental group of the Lagrangian orbit is done by several steps.\\
{\it Step 1}. The representation $\rho$ is of quaternionic type, so it preserves
a quaternionic structure $J \in \mbox{End}(\C^{32})$,
such that $J^2=-\mbox{id}$.
Denote by $\lambda$ the maximal weight of $\rho$ and by $T$ a fixed maximal
torus of $K$.
Note first that the Weyl group $W_{\Spin(12)}$ contains $-1$. Let $w\in N_{\Spin(12)}(T)$
induce $-1$ on $\mathfrak{t}$. Since $w(\lambda)=-\lambda$, $w$ preserves also $\lambda^\perp$,
therefore $w\in N_K(\SU(6))\subset N_K(\U(6))$.  On the other hand $w$ cannot lie in $\U(6)$
because otherwise $w$ should belong to $W_{\SU(6)}$ but ${-1}\notin W_{\SU(6)}$.
Hence $w$ generates $N(\U(6)/\U(6)$, and by \cite{BR} we
know that $N_{\SO(12)}(\U(6)/\U(6)\simeq \Z_2$. \\
{\it Step 2}. Take $p\in \P^{31}$ with $K\cdot p$ Lagrangian and $K_p^o$ locally
isomorphic to $\SU(6)$. Since $\rho$ is self dual, Remark \ref{selfdual} implies that
$p=u_1+u_{-1}$ with $u_{\pm 1} \in V_{\pm \lambda}$ and $\|u_1\|=\|u_{-1}\|$.
Now $K_p\subset N_{\SO(12)}(\SU(6))\subset N_{\SO(12)}(\U(6))$, hence if $k\in K_p$ then
$k\in w^i\U(6)=w^i T^1\cdot \SU(6)$ for $i=0,1$, where $T'$ is the center of $\U(6)$.\\
{\it Step 3}. Let $v_1\in V_{\lambda}$ be fixed  and take $v_2=Jv_1\in V_{-\lambda}$.
It is possible to choose $x\in\U(1)$ such that $p=v_1+xv_2$ and $w\cdot p\in \C\cdot p$, and the
$K$-orbit through $p$ is Lagrangian.\\
{\it Step 4}. Let $T^1\in\U(6)$ be the center of $\U(6)$. We consider the homomorphism
$c:T^1\to \U(1)$ such that, for every $t\in T^1$, $t\cdot v_1=c(t)\cdot v_1$,
with $c(t)\neq 1$. By easy computations we get that $t\cdot v_2=\overline{c(t)} v_2$.
Let $k\in K_p\subset w^i \cdot T^1\cdot \SU(6)=T^1\cdot w^i \SU(6)$, since both $w$ and
$\SU(6)$ fix $[p]$, then $k\in H\cdot w^i \cdot \SU(6)$ where $H:=\{t\in T^1, t[p]=[p]\}$.
Now $t\cdot p\in \C\cdot p$ if and only if $c(t)=\pm 1$ i.e.
$H=\Ker (c^2)\subset T^1$ is cyclic.\\
{\it Step 5}. Recall that $w^2\in T$, hence it commutes with $H$. In $\SO(12)$, $w$ can be taken as $\diag(B,B,\ldots,B)$ where $B=
\left[
\begin{array}{cc}
0 & 1 \\
1 & 0
\end{array}
\right]$; then $w\in \Spin (12)$, taken in $\pi^{-1}(w_{\SO(12))}$, is such that
$w^2\in \pi^{-1}(e)\simeq \Z_2$ i.e. $w^4=id$.\\
{\it Step 6}. Now we determine $H$. If $uJ$, with $u\in i\R$,
is a generic element of $\mathfrak{t}^1$, and
$J:=\diag(A,A,\ldots,A)$ where $A=
\left[
\begin{array}{cc}
0 & 1 \\
-1 & 0
\end{array}
\right]$
we have $H=\{uJ| \exp (uJ)v_1=\pm v_1\}$, and
recalling that  $\lambda=\frac{1}{2}(\omega_1+\omega_2+\cdots+\omega_6)$,
where $\omega_i$ are the fundamental weights, we get
$$H=\{uJ|u=\frac{\pi}{3}i\cdot k, k\in\Z\}.$$ Obviously $H\cap \SU(6)=\Z_6$ therefore
$[k]\in K_p/K_p^o$ is generated by $[w]$ and by $\alpha$, where $\alpha$ is a non trivial
 element of $H/H\cap \SU(6)\simeq \Z_2$. Now $w^2\in H$, thus $[w]^2$ equals $\alpha$ or $1$
and $|K_p/K_p^o|\leq 4$. The claim follows from the fact that the invariant has degree $d=4$.
\item {\it $K=\Spin(14)$, $\rho=\Lambda_e$}
Let $p$ be such that $K\cdot p$ is the Lagrangian $K$-orbit in
$\P(V)$. From \cite{SK} we know that $G_p^o$ is $\G_2^\C\times\G_2^\C$ and
from the inclusion
\[
\G_2\times\G_2 \subset \SO(7)\times\SO(7) \subset SO(14)
\]
which lifts to $\Spin(14)$ ($\G_2\times\G_2$ is simply connected),
we get $K_p^o=\G_2\times\G_2$. Now we claim that
$\pi_1(K_p)=K_p/K_p^o$ is exactly $\Z_8$. Since in this case the
degree of the invariant is 8 (see \cite{SK}) to prove this fact it
is sufficient to show that $|N_K(K_p^o)/K_p^o|$ is at most 8. \\
First we compute $N_{\SO(14)}(K_p^o)/K_p^o$. As an
automorphism of $\G_2\times\G_2$ an element $g$ of $N_{\SO(14)}(K_p^o)$ can
either preserve or interchange the $\G_2$ factors. Since
$\mathfrak{g}_2$ has no outer automorphism and the centralizer of
$\G_2$ in $\O(7)$ is $\{\pm \mbox{Id}\}$ we have that $N_{\SO(14)}(K_p^o)$
is given by the following four connected components:
\[
\left \{ \left[
\begin{array}{cc}
A & 0 \\
0 & B
\end{array}
\right]: A,B\in \G_2 \right\} \cup
\left \{ \left[
\begin{array}{cc}
-A & 0 \\
0 & -B
\end{array}
\right]: A,B\in\G_2 \right\} \cup
\]
\[
\left \{ \left[
\begin{array}{cc}
0 & A \\
-B & 0
\end{array}
\right]: A,B\in \G_2 \right \}\cup
\left
\{ \left[
\begin{array}{cc}
0 & -A \\
B & 0
\end{array}
\right]: A,B\in\G_2 \right \} .
\]
Thus $N_{\SO(14)}(K_p^o)/K_p^o \cong \Z_4$. Now note that the Lie group covering map
$\Spin(14)\to \SO(14)$ induces an epimorphism
\[
\frac{N_{\Spin(14)}(K_p^o)}{K_p^o} \to \frac{N_{\SO(14)}(K_p^o)}{K_p^o}
\]
whose kernel is $\Z_2$ and the claim follows.
\item {\it $K=\Ea_6$, $\rho=\Lambda_1$}. As before let $p$ be such that $K\cdot p$ is the Lagrangian $K$-orbit in
$\P(V)$. From \cite{SK} we know that $G_p^o$ is $\Fa_4^\C$ but $\Fa_4 \subset \Ea_6$, hence
$K_p^o=\Fa_4$. Following the same argument as in (\ref{SU(7)}), since $\Fa_4$ has only
inner automorphisms,
we get that $N_K (\Fa_4)$ is contained in $\Fa_4\cdot C_K(\Fa_4)$. Now $\Fa_4$ acts on
$\C^{27}=\C\oplus\C^{26}$ irreducibly on the second summand, hence $C_K(\Fa_4)$ acts on each summand
as scalar multiplication. Therefore $C_K(\Fa_4)$ is contained in a 2-dimensional torus
and $N_K(\Fa_4)/\Fa_4$ is abelian.
Since in this case the invariant has degree 3, we have $K_p=\Fa_4\cdot\Z_3$. Note that $\Z_3$ is the center of $\Ea_6$ which  acts trivially on $\
P^{26}$.
\item  {\it $K=\Ea_7$, $\rho=\Lambda_1$}. As before let $p$ be such that
$K\cdot p$ is the Lagrangian $K$-orbit in $\P(V)$. From \cite{SK} we know that
$G_p^o$ is $\Ea_6^\C$ but $\Ea_6 \subset \Ea_7$, hence $K_p^o=\Ea_6$.
This representation is self-dual.
\item{\it $K=\G_2$, $\rho=\Lambda_2$}. The orbits of $\G_2$ are the same of $\SO(7)$
(see case (\ref{Spin(7)})), therefore the Lagrangian orbit is
\[
\frac{\G_2}{\SU(3)\cdot\Z_2}=\frac{\SO(7)}{\SO(6)\cdot\Z_2}=\R\P^6.
\]


\end{enumerate}

We have thus proved Theorem \ref{classificazione}

\begin{center} { {\bf
Table 1:} Lagrangian orbits of simple Lie Groups in Projective spaces }
\end{center}
$$
\begin{array}{|l|l|l|l|l|l|l|l|l|l} \hline
 & \quad K &\quad \rho  &\dim_\C \P(V)&\text{cond.}& K_p^0& K_p/K_p^0& d \\\hline
1&\SU (n)    &  2 \Lambda_1          &  \  \frac{n(n+1)}{2}-1  &          & \SO(n)&  \Z_{n}             & n\\
2&\SU (n)   &  \Lambda_1\oplus \Lambda_1^*           &  \ 2n-1         &          & \SU(n-1) & \Z_2     &  2  \\
3&\SU (n)   &  \underbrace{\Lambda_1\oplus\cdots\oplus \Lambda_1}_n &\ n^2-1            &      &\{1\}    & \Z_n &   n  \\
4&\SU (2n)   &  \Lambda_2                     &  \ n(2n-1)-1   & n \geq 3 &  \Sp(n)&\Z_{2n}           & 2n \\
5&\SU (2n+1)&  \Lambda_2\oplus \Lambda_1      &\ 2n^2+3n+1     & n\geq 2  & \Sp(n) & \Z_{n+1} &  n+1 \\
6&\SU (2)    &  3 \Lambda_1                   &  \ 3           &          & \{1\}&\Z_3\rtimes \Z_4& 4\\
7&\SU (6)    &  \Lambda_3                     &  \ 19          &          & (\SU(3)\times\SU(3)) & \Z_4 & 4\\
8&\SU (7)    &  \Lambda_3                     &  \ 34          &          & \G_2 & \Z_7                 & 7\\
9&\SU (8)    &  \Lambda_3                     &  \ 55          &          & {\rm Ad}(\SU(3)) & \Z_{16}            & 16 \\
10&\Sp (n)    &  \Lambda_1\oplus \Lambda_1            & \ 4n-1            &          &\Sp(n-1) & \Z_2     &  2  \\
11&\Sp (3)    &  \Lambda_3                     &  \ 13          &         &\SU(3) & \Z_4                & 4  \\
12&\SO (n)    &  \Lambda_1                     &   \ n-1        &n \geq 3 & \SO(n-1) &  \Z_2            & 2\\
13&\Spin  (7) & {\rm spin \ rep.}               & \ 7            &        & \G_2 &  \Z_2                & 2\\
14&\Spin (9)  & {\rm spin \ rep.}              &  \ 15          &         &\Spin(7) &  \Z_2             & 2\\
15&\Spin (10) & \Lambda_e\oplus \Lambda_e      & \ 31              &      & \G_2 & -           &   4  \\
16&\Spin (11) & {\rm spin \ rep.}              &  \ 31          &        &  \SU(5) &  \Z_4                   & 4\\
17&\Spin (12) & \Lambda_e                      &  \ 31           &        &  \SU(6) &  \Z_4                  &  4\\
18&\Spin (14) & \Lambda_e                      & \ 63           &         & (\G_2\times\G_2) &  \Z_8    &  8  \\
19&\Ea_6      & \Lambda_1                     & \ 26           &          &  \Fa_4 & \Z_3               &  3 \\
20&\Ea_7      & \Lambda_1                      & \ 55           &         &     \Ea_6 &   -             & 4 \\
21&\G_2       & \Lambda_2                       & \ 6            &        &   \SU(3) &  \Z_2            &  2 \\
\hline
\end{array}
$$
In the Table the connected components of the isotropy subgroups
$K_p$ of points $p$ through which the $K$-orbit is Lagrangian are
listed in the fifth column.

\end{document}